\documentclass[11pt,a4paper]{article}
\usepackage[margin=27mm]{geometry}
\usepackage{amsmath,amssymb,amsthm,booktabs,thmtools}
\usepackage{microtype,needspace}
\usepackage{enumitem}
\usepackage[hidelinks]{hyperref}
\usepackage[nameinlink,noabbrev]{cleveref}
\newtheorem{theorem}{Theorem}
\newtheorem{conjecture}[theorem]{Conjecture}
\newtheorem{lemma}[theorem]{Lemma}

\theoremstyle{remark}

\crefname{theorem}{theorem}{theorems}
\crefname{conjecture}{conjecture}{conjectures}
\Crefname{conjecture}{Conjecture}{Conjectures}
\Crefname{theorem}{Theorem}{Theorems}
\crefname{lemma}{lemma}{lemmas}
\Crefname{lemma}{Lemma}{Lemmas}
\crefname{corollary}{corollary}{corollaries}
\Crefname{corollary}{Corollary}{Corollaries}
\crefname{remark}{remark}{remarks}
\Crefname{remark}{Remark}{Remarks}
\newlist{statementparts}{enumerate}{1}
\setlist[statementparts]{label=\textup{(\roman*)},ref=\thetheorem(\roman*)}
\crefname{statementpartsi}{lemma}{lemmas}
\Crefname{statementpartsi}{Lemma}{Lemmas}

\DeclareMathOperator{\rank}{rank}
\allowdisplaybreaks
\title{Rank and Order of Graphs: A Seidel Matrix Approach}
\author{Ebrahim Ghorbani\\[5pt]
\normalsize\sl Institute for Algorithms and Complexity\\
\normalsize\sl Hamburg University of Technology, Hamburg, Germany\\[2pt]
\small \texttt{ebrahim.ghorbani@tuhh.de}}
\date{}
\begin{document}
\maketitle
\begin{abstract}
Motivated by the problem of relating the chromatic
number of a graph to its rank and its connections
to communication complexity, Kotlov
and Lov\'asz~(1996) proved that the maximum order
of a reduced graph of rank $r$ is $\Theta(2^{r/2})$,
where the rank of a graph is the rank of its adjacency
matrix $A$, and reduced means that $A$ has distinct
nonzero rows.
Akbari, Cameron and Khosrovshahi~(2004) conjectured
that this maximum is $2^{r/2+1}-2$ for even $r$
and $5\cdot2^{(r-3)/2}-2$ for odd $r$.
Ghorbani, Mohammadian and Tayfeh-Rezaie~(2015)
proved that if the conjecture is false, then there
exists a counterexample of rank at most~$46$.
We prove the conjecture and characterize all graphs
of maximum order. Our proof is independent of the
reduction to rank at most $46$.
We also prove the analogous conjecture for the corank
of a graph, defined as the rank of $A+I$.
Our approach brings both conjectures into a unified
framework by recasting them as extremal problems for
Seidel matrices.
Our proofs combine linear algebra over the binary
and real fields with elementary discrete Fourier analysis.
\end{abstract}

{\small
	
	\par}

\section{Introduction}

Throughout, all graphs are finite and simple.
Let $A(G)$ denote the adjacency matrix of a graph~$G$.
The \emph{rank} of $G$ is defined as $\rank A(G)$,
where rank is taken over the reals unless otherwise specified.
A graph is called \emph{reduced} if its adjacency matrix
has distinct nonzero rows.

Kotlov and Lov\'asz \cite{KL} proved that the maximum
order of a reduced graph of rank $r$ grows as
$\Theta(2^{r/2})$. For each $r\ge2$, they constructed
 reduced graphs of rank $r$ and order
\[
 m(r)=
 \begin{cases}
  2^{r/2+1}-2,&r\text{ even},\\
  5\cdot2^{(r-3)/2}-2,&r\text{ odd}.
 \end{cases}
\]
We also put $m(1)=0$.
More generally, reduced graphs of rank $r$ and order
$m(r)$ can be obtained by successive applications of
the following two constructions, starting from $K_2$
for even $r$ and from $K_3$ for odd $r$.
For a graph $F$ with adjacency matrix $A=A(F)$, define
the graphs $F^{\rm a}$ and $F^{\rm c}$ by Construction
{\rm(a)} (see \cite{KL}) and Construction {\rm(c)}
(see \cite{ACK,HP}), respectively:
\begin{equation}\label{eq:constructions}
 A(F^{\rm a})=
 \begin{pmatrix}A&A&0&0\\A&A&\mathbf1&0\\0&\mathbf1^\top&0&1\\0&0&1&0\end{pmatrix},
 \quad\text{and}\quad
 A(F^{\rm c})=
 \begin{pmatrix}A&J-A&\mathbf1&0\\J-A&A&0&\mathbf1\\\mathbf1^\top&0&0&1\\0&\mathbf1^\top&1&0\end{pmatrix}.
\end{equation}
Here $J$ is the all-ones matrix of the same order as $A$,
and $\mathbf1$ denotes the all-ones column vector of the
appropriate length. 

\begin{conjecture}[{Akbari,	Cameron and Khosrovshahi~\cite{ACK}}]\label{conj:rank}
Let $G$ be a reduced graph of order $n$ and rank $r\ge2$. Then
$n\le m(r).$
Equality holds if and only if $G$ is obtained from
$K_2$ when $r$ is even, or from $K_3$ when $r$ is odd,
by successive applications of Constructions {\rm(a)} and {\rm(c)}.
In particular, up to isomorphism, there are exactly
$\lfloor r/2\rfloor$ graphs attaining equality.
\end{conjecture}

We have restated the conjecture in our terminology, since
the definition of a reduced graph in \cite{ACK} allows
one isolated vertex.
Haemers and Peeters \cite{HP} proved the conjecture
for graphs containing an induced matching of size $r/2$
when the rank $r$ is even, and containing an induced
subgraph consisting of the disjoint union of a matching
of size $(r-3)/2$ and a triangle when $r$ is odd.
Ghorbani, Mohammadian and Tayfeh-Rezaie
\cite{GMT} proved that if the inequality in
\Cref{conj:rank} fails, then there exists a
counterexample of rank at most $46$.

The \emph{corank} of $G$ is defined as
$\rank(A(G)+I)$, where $I$ is the identity matrix.
A graph is called \emph{coreduced} if $A(G)+I$
has distinct rows.
The corresponding problem for the corank of a graph
leads to the following conjecture.

\begin{conjecture}[{\cite{TR}, see also \cite{EGcorank}}]\label{conj:corank}
Let $G$ be a coreduced graph of order $n$ and corank
$r\ge2$. Then $n\le m(r)$.
For $r\ge4$, there is a unique graph attaining equality,
up to isomorphism, namely the complement of the graph
obtained from $K_2$ when $r$ is even, or from $K_3$
when $r$ is odd, by repeated applications of
Construction {\rm(c)}.
For $r=2$, equality holds only for $2K_1$, and for
$r=3$, only for $3K_1$ and $P_3$.
\end{conjecture}


The main result of this paper is the resolution of
\Cref{conj:rank,conj:corank}. In fact, we show that both conjectures can be treated
within a unified framework based on Seidel matrices,
and establish more general results from which both
conjectures follow. Our proofs are
entirely analytic and do not rely on the reduction
to rank at most $46$ established in~\cite{GMT}.

\subsection{Motivation and related work}

The problem of bounding the order of reduced graphs in terms
of their rank arose from the study of the relationship between
the chromatic number and adjacency rank. Van~Nuffelen~\cite{VN} conjectured
that $\chi(G)\le \rank A(G)$ for every graph $G$ with at
least one edge. After Alon and Seymour \cite{AS} disproved
this conjecture, the question became how large the gap
between these two parameters could be. Razborov \cite{Raz}
showed that the gap can be superlinear, and Nisan and
Wigderson \cite{NW} obtained a stronger separation.
It was in this context that Kotlov and Lov\'asz initiated
the study of the maximum order of reduced graphs of given
rank.

The problem is closely related to the log-rank conjecture
of Lov\'asz and Saks \cite{LS}, a fundamental open
problem in communication complexity. In its
graph-theoretic formulation, the conjecture asserts that
there exists an absolute constant $c$ such that
$ \log\chi(G)=O((\log\rank A(G))^c)$
for every graph $G$ with at least one edge; see
\cite{LS2}. Identifying vertices with identical
neighbourhoods and deleting isolated vertices preserves
both the rank and the chromatic number of such a graph.
It therefore suffices to study reduced graphs, which
provides a natural motivation for bounding their order
in terms of their rank.

For a Boolean communication matrix of rank $r$,
Lovett \cite{Lovett} established an upper bound of
$O(\sqrt{r}\log r)$ on the deterministic communication
complexity. Sudakov and Tomon \cite{ST} improved this
bound to $O(\sqrt{r})$ using discrepancy and spectral
methods. These results also bound $\log\chi(G)$:
a deterministic protocol for the adjacency matrix
partitions it into monochromatic rectangles, and the
vertices corresponding to the diagonal entries of
each rectangle form an independent set in $G$. On the lower-bound side,
G\"o\"os, Pitassi and Watson \cite{GPW} obtained a further
separation between deterministic communication
complexity and the logarithm of the rank.

More recently, Beniamini, Linial and Shraibman
\cite{BLS} investigated connections between the
log-rank conjecture and Ramsey theory through the
study of graphs with small clique number and small
rank of the adjacency matrix of their complement.
This provides another motivation for understanding
graphs of low rank.

The problem of determining the maximum order of reduced
graphs of a given rank has also been studied for restricted
graph classes. The maximum order and the corresponding
extremal graphs are known for trees \cite{GMTbip,EGHK} and bipartite graphs
\cite{GMTbip}, as well as for triangle-free graphs
\cite{GMTtri}. A closely related problem is the
classification of reduced graphs that are maximal
with respect to taking induced supergraphs without
increasing the rank; see \cite{EGHK}.

The rank and corank problems also arise naturally in
the study of eigenvalue multiplicities of graphs.
Given an eigenvalue $\mu$ of a graph $G$, a classical
problem in spectral graph theory is to bound the order
of $G$ in terms of $r=\rank(A(G)-\mu I)$, that is,
the codimension of the eigenspace corresponding to
$\mu$. For $\mu\notin\{0,-1\}$, Rowlinson \cite{Row}
obtained a quadratic upper bound on the order of $G$,
which was subsequently improved by Bell and Rowlinson
\cite{BR} to
$|V(G)|\le r(r+1)/2,$
apart from the exceptional cases where $\mu=1$ and
$G$ is $K_2$ or $2K_2$.
The eigenvalues $0$ and $-1$ are fundamentally different
in this respect. By duplicating vertices with their open
or closed neighbourhoods, respectively, one can increase
the order of a graph arbitrarily without increasing
its rank or corank. Excluding these operations leads
naturally to the rank problem for reduced graphs and
the corank problem for coreduced graphs. In contrast
to the quadratic bound for the other eigenvalues,
the optimal bounds in these two cases are exponential
in $r$.

Esmailian and Ghorbani \cite{EGLO} also established
connections between these two problems and a discrete
variation of the Littlewood--Offord problem concerning
subsums of vectors that lie in the binary cube.

\subsection{A unified Seidel matrix approach}

The \emph{Seidel matrix} of a graph $G$ is
$S(G)=J-I-2A(G)$.
We call $$U(G)=S(G)+I=J-2A(G)$$ the \emph{shifted Seidel matrix} of $G$.
Thus a shifted Seidel matrix is a symmetric matrix
with entries in $\{1,-1\}$ and all diagonal entries
equal to $1$.

A transformation $U\mapsto DUD$, where $D$ is a
diagonal matrix with diagonal entries in $\{1,-1\}$,
is called a \emph{sign congruence}.
Two shifted Seidel matrices are \emph{switching equivalent}
if one can be obtained from the other by a sign
congruence and a simultaneous permutation of rows
and columns.

For a shifted Seidel matrix $U$, partition its rows into
classes, where two rows belong to the same class if they
are equal up to sign. Let $q(U)$ denote the number of
these classes. A \emph{representative matrix} of $U$
is a principal submatrix of order $q(U)$ obtained by
retaining one index from each class.
We call $U$ \emph{reduced} if its rows are distinct
up to sign, or equivalently, if $q(U)$ equals the
order of $U$.

For a shifted Seidel matrix $P$, define
\begin{equation}\label{eq:phi}
	\Phi(P)=
	\begin{pmatrix}
		P&P&\mathbf1\\
		P&P&-\mathbf1\\
		\mathbf1^\top&-\mathbf1^\top&1
	\end{pmatrix}.
\end{equation}
Define $E_3=2I_3-J_3$, $E_4=2I_4-J_4$, and $E_r=\Phi(E_{r-2})$ for $r\ge5$.

\medskip
The following theorem provides an equivalent formulation
of \Cref{conj:rank} in terms of shifted Seidel matrices
and also forms the basis of our proof of \Cref{conj:corank}.

\begin{theorem}\label{thm:seidel}
	Let $U$ be a shifted Seidel matrix of rank $r+1$,
	where $r\ge2$. Then
	\[
	q(U)\le m(r)+1.
	\]
	Equality holds if and only if a representative
	matrix of $U$ is switching equivalent to $E_{r+1}$.
\end{theorem}

The connection with the two conjectures is simple.
Let $\widehat G$ be the graph obtained from $G$ by
adding an isolated vertex. We have
\[
\rank U(\widehat G)=\rank A(G)+1,\quad\text{and}\quad
\rank(U(\overline G)+J)=\rank(A(G)+I).
\]
The first identity follows by taking the Schur complement
at the added index, and the second from
$U(\overline G)+J=2(A(G)+I)$.
Moreover, if $G$ is reduced, then $U(\widehat G)$
is reduced, and if $G$ is coreduced, then
$U(\overline G)$ has distinct rows.
These identities allow us to translate both conjectures
into the setting of shifted Seidel matrices.
While the order bound in \Cref{conj:rank} follows directly
from \Cref{thm:seidel}, the corank problem requires
an additional analysis of the rank-one perturbation $U+J$,
which we carry out in \Cref{sec:corank}.
The complete proofs of \Cref{conj:rank,conj:corank},
including the equality characterizations, are given in
\Cref{sec:complete-rank,sec:corank}, respectively.

\paragraph{Outline.}
The proof of \Cref{thm:seidel} proceeds by induction on the
rank of $U$. The key parameter is $\delta_*(U)$, the minimum
number of indices whose deletion makes a reduced shifted
Seidel matrix $U$ nonreduced.

We first encode row differences in a vector space over
$\mathbb F_2$ whose dimension is bounded in terms of the
real rank of $U$ (\Cref{lem:binary-space}).
Combining elementary discrete Fourier analysis with linear
algebra over the reals, we derive inequalities relating
the rank of $U$ to the distribution of its row differences
(\Cref{lem:fourth-moment,lem:general-moment}).
A key step is to bound the number of pairs of rows with
the same binary difference (\Cref{lem:multiplicity}).
Together, these results yield an upper bound on
$\delta_*(U)$ (\Cref{thm:tau-cases}).

We then study the principal submatrix obtained by deleting
a minimum set of indices that makes $U$ nonreduced.
We show that its rank drops by at least two and that
each of its row classes has at most two members
(\Cref{lem:minimum-deletion}).
Using rank inequalities and Schur complements, we establish
that $\delta_*(U)=1$ in the cases needed for the proof
of \Cref{thm:seidel} (\Cref{lem:small-tau}).
This gives both the order bound and the recursive structure
of the extremal matrices, from which \Cref{conj:rank} follows.

Finally, to prove \Cref{conj:corank}, we show that reduced
shifted Seidel matrices of sufficiently large order
do not admit certain rank-preserving integral extensions
(\Cref{lem:integral-q}).
Combining this result with \Cref{thm:seidel}, we obtain
a sharp bound on the order of $U$ in terms of the rank
of $U+J$ and characterize the equality cases
(\Cref{thm:perturbation}).


\section{Preliminaries}\label{sec:prelim}

For a matrix $U$ indexed by $Y$ and $R\subseteq Y$, let $U[R]$ be the principal submatrix indexed by $R$. We use bold lowercase letters for vectors and ordinary letters for their coordinates: $x_i$ denotes the coordinate of $\mathbf x$ indexed by $i$.

Every representative matrix of a shifted Seidel matrix
$U$ has the same rank as $U$, since the discarded rows
and columns are signed copies of the retained ones.
The representative matrix is reduced because two
retained rows that agree up to sign on the retained
columns must also agree up to sign on the deleted
columns.
Moreover, any two representative matrices of $U$ are
switching equivalent.

We use two standard rank identities. For a real symmetric block matrix $M=\left(\begin{smallmatrix}B&C\\C^\top&D\end{smallmatrix}\right)$, the Frobenius rank inequality gives \cite[0.4.5(e)]{HJ}
\begin{equation}\label{eq:row-rank}
 2\rank(B\ C)\le \rank M+\rank B.
\end{equation}
If $B$ is nonsingular, the matrix $D-C^\top B^{-1}C$ is called the \emph{Schur complement} of $B$ in $M$, and \cite[0.8.5(b)]{HJ}
\[
 \rank M=\rank B+\rank(D-C^\top B^{-1}C).
\]
We will also use the elementary identity
\begin{equation}\label{eq:J-I-inverse}
 (J_l-I_l)^{-1}=\frac{1}{l-1}J_l-I_l
 \quad\text{for }l\ge2.
\end{equation}
Every real symmetric matrix has a nonsingular principal
submatrix of order equal to its rank \cite[Section 0.7.6]{HJ}. If $U$ is reduced and $U[R]$ has the same rank as $U$, restriction of the row space to $R$ is injective. Hence $U[R]$ is reduced.

\subsection{Switching and deletion}

Let $U=(u_{ij})$ be indexed by $Y$. For $z\in Y$, let $D_z$ be the diagonal matrix obtained by placing
row $z$ of $U$ on the diagonal. The matrix $D_zUD_z$ has an all-ones row and column at $z$; we call this operation \emph{normalization} at $z$. The $(i,j)$-entry of the normalized matrix is $u_{zi}u_{ij}u_{zj}$, which is invariant under sign congruence. On putting $z$ first, the normalized matrix has the form
\begin{equation}\label{eq:normal-form}
 \begin{pmatrix}1&\mathbf1^\top\\\mathbf1&J-2C\end{pmatrix},
 \quad\text{with}\quad \rank U=1+\rank C,
\end{equation}
where $C$ is a symmetric binary matrix with zero diagonal; the rank identity follows by the Schur complement. The matrix $U$ is reduced if and only if the rows of $C$ are distinct and nonzero. If $U$ has rank one, $C$ is zero, so every row of the normalized $U$ is the all-ones row and $q(U)=1$. In particular, $U$ cannot have rank two, since any
nonzero symmetric binary matrix with zero diagonal
has a nonsingular principal submatrix of order two.

For distinct $i,j\in Y$, let $D_{ij}$ be the smallest
subset of $Y\setminus\{i,j\}$ whose deletion makes
rows $i$ and $j$ equal up to sign.
Equivalently,
\[
 D_{ij}=\{\ell\in Y:u_{ij}u_{i\ell}u_{j\ell}=-1\}.
\]
Put
\[
 \delta_*(U)=\min_{i\ne j}|D_{ij}|,
\]
where the minimum is used for matrices of order at least two.
Thus $\delta_*(U)$ is the minimum number of indices
whose deletion makes the resulting principal submatrix
nonreduced. The sets $D_{ij}$ are invariant under sign
congruence. After normalization at $i$, rows $i$ and $j$
differ precisely on $D_{ij}$. In particular, $U$ is
reduced if and only if $\delta_*(U)>0$.

\begin{lemma}\label{lem:deletion}
Let $U$ be reduced of rank $k$ and order at least two.
\begin{statementparts}
\item\label{lem:deletion:reduced} If $R\subseteq Y$ and $|Y\setminus R|<\delta_*(U)$,
then $U[R]$ is reduced.
\item\label{lem:deletion:rank} For distinct $i,j$, if $R=Y\setminus D_{ij}$, then $\rank U[R]\le k-2$.
\end{statementparts}
\end{lemma}
\begin{proof}
	The reducedness assertion follows directly from the definition of
	$\delta_*(U)$.
	
	For the rank bound, normalize $U$ at $i$. Since $U$ is reduced,
	rows $i$ and $j$ of $U$ are distinct, whereas the
	corresponding rows of $U[R]$ are equal.
	Thus, if $\ell$ denotes the rank of the rows of $U$
	indexed by $R$, then $\ell\ge\rank U[R]+1$.
	By \eqref{eq:row-rank}, we have
	$2(\rank U[R]+1)\le 2\ell\le k+\rank U[R],$
	which yields $\rank U[R]\le k-2$.
\end{proof}

\subsection{The construction $\Phi$}

For a shifted Seidel matrix $P$ of order $l$, we call the
indices $i$ and $l+i$ of $\Phi(P)$ a pair, for $1\le i\le l$.
\Needspace{8\baselineskip}
\begin{lemma}\label{lem:phi}
Let $P$ be a shifted Seidel matrix of order $l$ and rank $p$.
\begin{statementparts}
\item\label{lem:phi:rank} $\Phi(P)$ has order $2l+1$ and rank $p+2$.
\item\label{lem:phi:reduced} If $P$ is reduced, then $\Phi(P)$ is reduced.
\item\label{lem:phi:switching} If $P$ and $P'$ are switching equivalent, then so are
$\Phi(P)$ and $\Phi(P')$.
\end{statementparts}
\end{lemma}
\begin{proof}
\medskip
\noindent\textit{Part} (i).
The order of $\Phi(P)$ is clear from its definition.
On the subspace of vectors of the form $(\mathbf x,\mathbf x,0)$,
$\Phi(P)$ acts as $2P$ and hence has rank $p$.
Its restriction to the subspace spanned by
$(\mathbf1,-\mathbf1,0)$ and the last standard basis vector
has rank two, while it vanishes on the subspace of vectors
of the form $(\mathbf x,-\mathbf x,0)$ with
$\mathbf1^\top\mathbf x=0$. These three invariant subspaces
form a direct sum of the entire space, so $\rank\Phi(P)=p+2$.

\medskip
\noindent\textit{Part} (ii).
Suppose $P$ is reduced. Rows from distinct pairs
cannot agree up to sign, since their restrictions to the
first $l$ columns are distinct rows of $P$ up to sign.
The two rows in a pair agree except at the last index.
The last row takes opposite values on each pair of columns,
whereas all other rows take equal values on each pair.
Thus $\Phi(P)$ is reduced.

\medskip
\noindent\textit{Part} (iii).
A simultaneous sign change of a row and column
of $P$ can be effected in $\Phi(P)$ by changing the signs
of both corresponding rows and columns and then
interchanging the two corresponding indices. Permutations
of the indices of $P$ act simultaneously on the pairs.
Hence $\Phi(P)$ and $\Phi(P')$ are switching equivalent.
\end{proof}

Since $E_3$ and $E_4$ are reduced and nonsingular,
\Cref{lem:phi:rank,lem:phi:reduced} and the recurrence
$m(k-1)=2m(k-3)+2$ for $k\ge5$ show that $E_k$
is reduced of rank $k$ and order $m(k-1)+1$ for every $k\ge3$.

For a shifted Seidel matrix $P$, put
\[
 P^\pm=\begin{pmatrix}P&-P\\-P&P\end{pmatrix}.
\]
The matrix $P^\pm$ has twice the order of $P$
and the same rank as $P$.
If $P$ is reduced, then $P^\pm$ has distinct rows.
Moreover, if $U$ has distinct rows and every row
class has two members, then $U$ is equal to $P^\pm$
up to simultaneous permutation, where $P$ is a
representative matrix of $U$.
If $P$ and $P'$ are switching equivalent, then $P^\pm$ and
$(P')^\pm$ are equal up to simultaneous permutation, since
in $P^\pm$ a sign change in $P$ corresponds simply to
interchanging the two corresponding indices, and permutations
of $P$ act simultaneously on its two copies.

We will use one more property of $\Phi(P)$ when $l\ge2$:
\begin{equation}\label{eq:no-full-difference}
 |D_{ij}|<2l-1\quad \text{for }i\ne j.
\end{equation}
For two indices in the same pair, $D_{ij}$ consists of the last index only. For indices in different pairs, the partner of either index does not belong to $D_{ij}$. If one index is the last index and the other is in a pair, the two indices in any other pair give opposite signs of the product defining $D_{ij}$, so one of them does not belong to $D_{ij}$. These cases prove the strict inequality, since $2l-1$ is the number of indices other than $i,j$.

\section{Bounding $\delta_*$ via linear algebra over $\mathbb F_2$}
\label{sec:eq-high}

Our aim in this section is to obtain an upper bound
on $\delta_*(U)$ for reduced shifted Seidel matrices
of sufficiently large order relative to their rank.

We encode row differences as elements of a vector
space over $\mathbb F_2$ whose dimension is controlled
by the real rank of $U$. Combining this approach
with bounds on $\operatorname{tr}(U^4)$, we obtain
the desired upper bound on $\delta_*(U)$.

\subsection{A switching-invariant binary space}

In the matrix arguments below, we write $k$ for the rank of $U$, so the bound in \Cref{thm:seidel} is $m(k-1)+1$.
Let $U$ be reduced of order $N$, rank $k$, and index set $Y$. Write $\mathbf b_i$ for row $i$ of $(J-U)/2$, viewed over $\mathbb F_2$. Let $\mathcal Q=\mathbb F_2^Y/\langle\mathbf1\rangle$
be the quotient vector space, and denote by $[\mathbf x]$
the image of $\mathbf x\in\mathbb F_2^Y$ in $\mathcal Q$.
Let $\mathcal V$ be the subspace of $\mathcal Q$ spanned
by the vectors $[\mathbf b_i+\mathbf b_j]$ for $i,j\in Y$,
and put $s=\dim\mathcal V$.
Put
\[
 s_*:=2\left\lfloor\frac{k-1}{2}\right\rfloor.
\]
For $X=[\mathbf x]\in\mathcal V$, define
\[
 f(X)=\left|\sum_{i\in Y}(-1)^{x_i}\right|,
 \quad\text{and}\quad
 \mu(X)=\left|\{(i,j)\in Y^2:[\mathbf b_i+\mathbf b_j]=X\}\right|.
\]
The function $f$ is well-defined, since replacing
$\mathbf x$ by $\mathbf x+\mathbf1$ only changes
the sign of the sum inside the absolute value.

\begin{lemma}\label{lem:binary-space}
With the above notation, the following hold.
\begin{statementparts}
\item\label{lem:binary-space:invariance}
The vectors $[\mathbf b_i+\mathbf b_j]$, the space $\mathcal V$,
and the functions $f$ and $\mu$ are invariant under sign congruence.
\item\label{lem:binary-space:dimension}
After normalization at any index $z$, the quotient map restricts
to an isomorphism from the binary row space of $(J-U)/2$
onto $\mathcal V$. In particular, $s$ is even and $s\le s_*$.
\item\label{lem:binary-space:multiplicity}
The vectors $[\mathbf b_i]$, $i\in Y$, are distinct,
$ f(0)=\mu(0)=N$, $\sum_{X\in\mathcal V}\mu(X)=N^2$, and  $\mu(X)\le N$ for $X\in\mathcal V$.
\end{statementparts}
\end{lemma}
\begin{proof}
\medskip
\noindent\textit{Part} (i).
Let $\mathbf a=(a_i)_{i\in Y}$
be a binary vector. A sign congruence with diagonal entries
$(-1)^{a_i}$ changes $\mathbf b_i$ to
$\mathbf b_i+\mathbf a+a_i\mathbf1$ over $\mathbb F_2$.
Hence $[\mathbf b_i+\mathbf b_j]$ is unchanged, as are
$\mathcal V$ and $\mu$. Since $f$ depends only on $X$,
it is also invariant under sign congruence.

\medskip
\noindent\textit{Part} (ii).
After normalization at $z$, every binary row has a zero
at $z$, and the row at $z$ is zero. Consequently the quotient
map is injective on the binary row space, and the images
of the rows span $\mathcal V$. Thus $s$ is the binary rank
of the matrix $C$ in \eqref{eq:normal-form}.
Since $C$ is symmetric with zero diagonal, its binary rank
is even: over $\mathbb F_2$, taking the Schur complement
of a nonsingular principal block of order two preserves
symmetry and the zero diagonal and reduces the rank by two.
Every nonzero symmetric matrix with zero diagonal has such
a block, so the assertion follows by induction. Its binary rank is at most
its real rank, since any minor with odd determinant is
nonsingular over $\mathbb R$. By \eqref{eq:normal-form},
$\rank C=k-1$, so $s\le s_*$.

\medskip
\noindent\textit{Part} (iii).
Since $U$ is reduced, the vectors $[\mathbf b_i]$
are distinct. Hence only pairs $(i,i)$ contribute to $\mu(0)$,
and $f(0)=N$ directly from its definition.
Every ordered pair contributes to exactly one $\mu(X)$,
so their sum is $N^2$. For each fixed $i$ and $X$ there is
at most one $j$ with $[\mathbf b_i+\mathbf b_j]=X$,
which gives $\mu(X)\le N$.
\end{proof}

\begin{lemma}\label{lem:fourth-moment}
With the above notation,
\begin{equation}\label{eq:parseval}
 \sum_{X\in\mathcal V}f(X)^2=2^sN
\end{equation}
and
\begin{equation}\label{eq:fourth}
 \frac{N^4}{k}\le\operatorname{tr}(U^4)
 =\sum_{X\in\mathcal V}\mu(X)f(X)^2.
\end{equation}
In particular,
\begin{equation}\label{eq:binary-bound}
 N^2\le k2^s.
\end{equation}
\end{lemma}
\begin{proof}
Normalize at any index. By \Cref{lem:binary-space:dimension}, the images of the binary rows belong to $\mathcal V$. For $i\ne j$, the linear functional $[\mathbf x]\mapsto x_i+x_j$ on $\mathcal V$ is nonzero, since the distinct columns of the normalized binary matrix differ on some row. Hence it takes each value equally often, giving
\[
 \sum_{X\in\mathcal V}f(X)^2
 =\sum_{i,j\in Y}\sum_{[\mathbf x]\in\mathcal V}(-1)^{x_i+x_j}
 =2^sN.
\]
Since $U$ is symmetric, we have
\[
 |(U^2)_{ij}|=f([\mathbf b_i+\mathbf b_j]).
\]
Consequently,
\[
 \operatorname{tr}(U^4)
 =\sum_{i,j\in Y}(U^2)_{ij}^2
 =\sum_{X\in\mathcal V}\mu(X)f(X)^2.
\]

Moreover, $\operatorname{tr}(U^2)=N^2$ and $\rank U=k$,
so Cauchy--Schwarz applied to the squares of the nonzero
eigenvalues gives the inequality in \eqref{eq:fourth}.
Finally, the bound $\mu(X)\le N$ in
\Cref{lem:binary-space:multiplicity} yields
\eqref{eq:binary-bound}.
\end{proof}

\begin{lemma}\label{lem:general-moment}
Let $h,L\ge0$. If $\mu(X)\le L$ for all nonzero $X\in\mathcal V$ with $f(X)>h$, then
\begin{equation}\label{eq:general-moment}
 \frac{N^4}{k}\le N^3+h^2N(N-1)+L(2^{s_*}N-N^2).
\end{equation}
\end{lemma}
\begin{proof}
	By \eqref{eq:fourth}, we have
	\[
	\frac{N^4}{k}
	\le N^3
	+\sum_{\substack{X\ne0\\ f(X)\le h}}\mu(X)f(X)^2
	+\sum_{\substack{X\ne0\\ f(X)>h}}\mu(X)f(X)^2.
	\]
	By \Cref{lem:binary-space:multiplicity}, the first sum
	is at most
	$h^2\sum_{X\ne0}\mu(X)=h^2N(N-1).$
	By the assumption on $\mu(X)$ and \eqref{eq:parseval},
	the second sum is at most
	$L\sum_{X\ne0}f(X)^2=L(2^sN-N^2).$
	Combining these bounds and using $s\le s_*$ from
	\Cref{lem:binary-space:dimension} yields the desired
	inequality.
\end{proof}

\subsection{Counting pairs with the same row difference}

The inequality in \Cref{lem:general-moment} becomes
stronger if we can bound $\mu(X)$ for vectors $X$
with large $f(X)$. Recall that a large value of $f(X)$
means that the corresponding rows differ in few
positions up to sign. The following lemma provides
the required bound on $\mu(X)$.

\begin{lemma}\label{lem:multiplicity}
Let $k\ge5$ and assume that the inequality in
\Cref{thm:seidel} holds for all ranks less than $k$.
Let $X\in\mathcal V$ be nonzero with $\mu(X)>0$, and put
$d=(N-f(X))/2$. Suppose
\[
 d<\delta_*(U)
 \quad\text{and}\quad
 N-d>m(k-3)+1.
\]
Then the following hold.
\begin{statementparts}
\item\label{lem:multiplicity:submatrix}
$U$ has a reduced principal submatrix of order
$N-d$, rank $k-1$, and binary space of dimension
at most $s-2$.
\item\label{lem:multiplicity:support}
$d\ge N-m(k-2)-1$.
\item\label{lem:multiplicity:pairs}
If $d>N-2\delta_*(U)$, then $\mu(X)=2$.
\end{statementparts}
\end{lemma}
\begin{proof}
\medskip
\noindent\textit{Parts {\rm(i)} and {\rm(ii)}}.
Choose a representative of $X$ with support $D$ of size
$d$ and put $R=Y\setminus D$.
We call a pair $(i,j)$ realizing if
$[\mathbf b_i+\mathbf b_j]=X$.
Normalize at an index of $R$. By \Cref{lem:binary-space:invariance}, the realizing pairs are unchanged. The binary row sum of every such pair is now precisely this representative, since it vanishes at the normalization index. Since $d<\delta_*(U)$, \Cref{lem:deletion:reduced}
implies that $U[R]$ is reduced.

Every realizing pair $i,j$ has both indices in $D$. Indeed, the binary row sum takes the same value at $i$ and $j$ by symmetry and the zero diagonal of the binary matrix. If both were in $R$, the corresponding rows of $U[R]$ would be equal. Let $\mathbf u_i$ denote row $i$ of $U$. The row difference $\mathbf u_i-\mathbf u_j$ is nonzero and supported on $D$. Consequently restriction of the row space of $U$ to $R$ has dimension at most $k-1$, giving $\rank U[R]\le k-1$. Since $N-d>m(k-3)+1$, the induction hypothesis gives
$\rank U[R]=k-1$ and $N-d\le m(k-2)+1$.
The latter inequality gives the claimed lower bound on $d$.

Under the restriction map
\[
 \mathbb F_2^Y/\langle\mathbf1\rangle
 \longrightarrow
 \mathbb F_2^R/\langle\mathbf1\rangle,
\]
the image of $\mathcal V$ contains the binary space
of $U[R]$. Since $X\ne0$ belongs to the kernel,
this binary space has dimension at most $s-1$.
By \Cref{lem:binary-space:dimension}, its dimension
is even, and hence at most $s-2$.

\medskip
\noindent\textit{Part} (iii).
Since $\rank U[R]=k-1$, the subspace of the row
space of $U$ consisting of vectors supported on $D$
has dimension one. All real row differences of realizing pairs are therefore equal up to sign (their nonzero entries are $\pm2$). If $\mu(X)>2$, then by \Cref{lem:binary-space:multiplicity}
there are two disjoint realizing pairs.
By interchanging their indices if necessary, we may
write them as $(i,j)$ and $(p,q)$ with
\[
 \mathbf u_i-\mathbf u_j=\mathbf u_p-\mathbf u_q.
\]
On $D$ the rows of $i,p$ agree and those of $j,q$ agree with the opposite values. All four indices lie in $D$, so $u_{ip}=1$ and $u_{iq}=-1$. On $R$, the rows in each realizing pair agree. If the rows of $i,p$ differ in $t$ coordinates, the rows of $i,q$ differ in $t+d$ coordinates. Hence
\[
 \delta_*(U)\le |D_{ip}|=t,
 \quad\text{and}\quad
 \delta_*(U)\le |D_{iq}|=N-t-d.
\]
Adding these inequalities gives $d\le N-2\delta_*(U)$.
Thus if $d>N-2\delta_*(U)$, then $\mu(X)\le2$.
Since $\mu(X)>0$ and each realizing pair occurs in both
orders, we have $\mu(X)=2$.
\end{proof}

We are now ready to combine the preceding results
to establish the desired upper bound on $\delta_*(U)$.

\begin{lemma}\label{thm:tau-cases}
Let $k\ge5$ and assume that the inequality in
\Cref{thm:seidel} holds for matrices of rank less than $k$.
If $U$ is reduced of rank $k$ and order
$N\in\bigl\{m(k-1),m(k-1)+1,m(k-1)+2\bigr\}$, then
\[
 \delta_*(U)\le m(k-3).
\]
\end{lemma}
\begin{proof}
	Put $B=m(k-3)+1$ and suppose, to the contrary,
	that $\delta_*(U)\ge B$. Note that $m(k-1)=2B$.
	
	\medskip
	\noindent\textit{Case} 1: $k=5$ and $N=6$.
	Normalize $U$
	at any index and let $C$ be the matrix in
	\eqref{eq:normal-form}. Since $\delta_*(U)\ge3$,
	each row of $C$ has at least three ones.
	If $c_{ij}=0$ for some distinct $i,j$, then rows $i$
	and $j$ must have ones at all three other indices.
	These two rows are therefore equal, contradicting
	the reducedness of $U$. Hence $C=J-I$ has rank five,
	contrary to $\rank C=k-1=4$.
	
	For Cases 2 and 3, put $h=N-2B\in\{0,1,2\}$.
	Let $X\in\mathcal V$ be nonzero with $\mu(X)>0$
	and $f(X)>h$, and put $d=(N-f(X))/2$. Then
	\[
	d<B\le\delta_*(U),\quad\text{and}\quad
	N-d>N-B\ge m(k-3)+1.
	\]
	Thus \Cref{lem:multiplicity} applies.
	
	\medskip
	\noindent\textit{Case} 2: $k\ge7$.
	By
	\Cref{lem:multiplicity:support},
	\[
	d-h\ge N-m(k-2)-1-h
	=m(k-1)-m(k-2)-1>0.
	\]
	Since $h\ge N-2\delta_*(U)$,
	\Cref{lem:multiplicity:pairs} gives $\mu(X)=2$.
	We may therefore apply \Cref{lem:general-moment}
	with $L=2$. Using $h\le2$ and
	\[
	2^{s_*}\le\frac{(m(k-1)+2)^2}{4}
	\le\frac{(N+2)^2}{4},
	\]
	we obtain
	\[
	\frac Nk
	\le 1+\frac{4(N-1)}{N^2}
	+\frac{2(2^{s_*}-N)}{N^2}
	\le\frac32+\frac4N-\frac2{N^2}<2.
	\]
	This contradicts $N\ge m(k-1)\ge2k$.
	
	\medskip
	\noindent\textit{Case} 3: $k\in\{5,6\}$ and $(k,N)\ne(5,6)$.
	We show that no nonzero $X\in\mathcal V$ satisfies
	both $\mu(X)>0$ and $f(X)>h$.
	
	If $k=5$, then $d<B=3$ and $N\ge7$.
	By \Cref{lem:multiplicity:submatrix}, $U$ has a reduced
	principal submatrix of rank four and order at least
	five, contradicting the induction hypothesis.
	
	If $k=6$, then $d<B=4$ and $N\ge8$.
	By \Cref{lem:multiplicity:submatrix}, $U$ has a reduced
	principal submatrix of order at least five whose
	binary space has dimension at most $s-2\le2$.
	By \Cref{lem:binary-space:dimension,lem:binary-space:multiplicity},
	after normalization its distinct binary rows lie in
	a vector space containing at most four vectors,
	a contradiction.
	
	Consequently, we may apply \Cref{lem:general-moment}
	with $L=0$, obtaining
	\[
	\frac Nk\le1+\frac{h^2(N-1)}{N^2}.
	\]
	The possible pairs $(N,h)$ are $(7,1)$ and $(8,2)$
	for $k=5$, and $(8,0)$, $(9,1)$, and $(10,2)$
	for $k=6$. If $h\le1$, the left-hand side is at least
	$4/3$, whereas the right-hand side is less than $8/7$.
	If $h=2$, the left-hand side is at least $8/5$,
	whereas the right-hand side is less than $3/2$.
	In either case we obtain a contradiction.
\end{proof}

\section{Extremal matrices and proofs of \Cref{thm:seidel} and \Cref{conj:rank}}
\label{sec:small}

Building on the upper bound on $\delta_*(U)$ established
in \Cref{sec:eq-high}, we determine the possible ranks
and row class sizes of the principal submatrix obtained
by deleting a set $D_{ij}$ of minimum size.
The structure of its row classes allows us to complete
the proof of \Cref{thm:seidel}. We then use this theorem
to prove \Cref{conj:rank}, including the characterization
of the extremal graphs.

\subsection{Reduction to $\delta_*(U)=1$}
\label{sec:delta-one}
The following lemma is related to
\cite[Lemma 4]{GMT}, but is formulated in terms of
shifted Seidel matrices.

\begin{lemma}\label{lem:minimum-deletion}
	Let $U$ be a reduced shifted Seidel matrix of rank $k\ge5$
	and order $N$, indexed by $Y$.
	Assume that the inequality in \Cref{thm:seidel} holds
	for matrices of rank less than $k$.
	Choose $i,j\in Y$ such that $D=D_{ij}$ has size
	$d=\delta_*(U)$, and put $R=Y\setminus D$.
	Suppose $N-d>m(k-3)+1$ and normalize $U$ at $i$.
	Then the following hold.
	\begin{statementparts}
		\item\label{lem:minimum-deletion:rank}
		$\rank U[R]\in\{k-3,k-2\}$.
		\item\label{lem:minimum-deletion:classes}
		Each row class of $U[R]$ has at most two members.
		If a class has two members, their corresponding rows
		in $U$ agree on $R$, while on $D$ one row consists
		entirely of ones and the other entirely of minus ones.
	\end{statementparts}
\end{lemma}
\begin{proof}
	\medskip
	\noindent\textit{Part} (i).
	By \Cref{lem:deletion:rank}, we have $\rank U[R]\le k-2$.
	Suppose $\rank U[R]\le k-4$ and choose an index $w\in D$.
	Then
	\[
	\rank U[R\cup\{w\}]\le\rank U[R]+2\le k-2.
	\]
	Since $|D\setminus\{w\}|<\delta_*(U)$,
	\Cref{lem:deletion:reduced} implies that
	$U[R\cup\{w\}]$ is reduced.
	Its order is greater than $m(k-3)+1$, contradicting
	the induction hypothesis.
	
	\medskip
	\noindent\textit{Part} (ii).
	Put $p=\rank U[R]$ and let $\ell$ be the rank of the
	rows of $U$ indexed by $R$.
	By \eqref{eq:row-rank}, we have $2\ell\le k+p$.
	Since $p\ge k-3$ and $\ell-p$ is an integer, it follows
	that $\ell-p\le1$.
	
	Let $\mathcal L$ be the space spanned by the rows of
	$U$ indexed by $R$. The restriction map from
	$\mathcal L$ to the row space of $U[R]$ has kernel
	of dimension $\ell-p$.
	Denote row $v$ of $U$ by $\mathbf u_v$.
	The nonzero vector $\mathbf u_i-\mathbf u_j$ belongs
	to this kernel: it is zero on $R$ and has every entry
	equal to two on $D$.
	Consequently, the kernel is one-dimensional and is
	spanned by $\mathbf u_i-\mathbf u_j$.
	
	Since column $i$ of $U[R]$ consists entirely of ones,
	rows in the same class of $U[R]$ are equal.
	Thus, if $v$ and $w$ are distinct indices belonging
	to the same class, then
	$\mathbf u_v-\mathbf u_w$ is a nonzero multiple of
	$\mathbf u_i-\mathbf u_j$.
	As the entries of $U$ are $\pm1$, one of these rows
	must be all ones on $D$ and the other all minus ones
	on $D$.
	A row class cannot contain three members, since two
	of their corresponding rows in $U$ would then agree
	on both $R$ and $D$, contradicting the reducedness
	of~$U$.
\end{proof}

By \Cref{lem:minimum-deletion}, deleting a smallest
set $D_{ij}$ leaves a principal submatrix whose row
classes have at most two members.
We now use the structure of these classes to show
that, under suitable conditions, deleting a single
index is sufficient to obtain a nonreduced matrix.
The only possible exception involves the matrix
\begin{equation}\label{eq:exception}
 W=
 \begin{pmatrix}
 1&\mathbf1^\top\\
 \mathbf1&J-2C_0
 \end{pmatrix},
 \quad\text{with}\quad
 C_0=\operatorname{diag}(J_2-I_2,J_3-I_3).
\end{equation}
By \eqref{eq:normal-form}, $W$ has order and rank six.

\begin{lemma}\label{lem:small-tau}
Suppose the hypotheses of \Cref{lem:minimum-deletion}
hold and that \Cref{thm:seidel} holds for matrices of
rank less than $k$. If
\[
 d\le m(k-3)
 \quad\text{and}\quad
 N-q(U[R])>m(k-3)+1,
\]
then $d=1$, unless $k=6$, $N=8$, and $U$ contains
a principal submatrix switching equivalent to~$W$.
\end{lemma}
\begin{proof}
	Put $B=m(k-3)+1$ and $q=q(U[R])$.
	By \Cref{lem:minimum-deletion:classes}, each row class
	of $U[R]$ has at most two members.
	Let $S$ contain the member of each two-element class
	whose row in $U$ is all ones on $D$, and let $S'$
	contain the other member.
	Let $T$ consist of the indices belonging to singleton
	classes. Put $a=|S|$ and $t=|T|$.
	Note that $i\in S$ and $j\in S'$.
	
	Since $\rank U[R]\le k-2$, the induction hypothesis
	gives
	\[
	q=a+t\le B,\quad\text{and}\quad N-d=2a+t.
	\]
	Together with $q\le N-B-1$, these relations yield
	\[
	a+d=N-q\ge B+1,\quad\text{and}\quad
	t=2q-N+d\le d-1.
	\]
	In particular, $a\ge2$.
	
	We first show that $U[S]$ and $U[D\cup\{i,j\}]$
	are reduced.
	Since both matrices have an all-ones column at $i$,
	it suffices to show that their rows are distinct.
	
	If two rows of $U[S]$ are equal, then the corresponding
	rows of $U$ also agree on $S'$, since each member of $S'$
	belongs to the same row class of $U[R]$ as its partner
	in $S$. They also agree on $D$ by the definition of $S$.
	Thus they can differ only on $T$, contradicting
	$t<\delta_*(U)$.
	
	In $U[D\cup\{i,j\}]$, rows $i$ and $j$ differ on $D$,
	and every row indexed by $D$ differs from both of
	these rows at column $j$.
	If two distinct rows indexed by $D$ agree on
	$D\cup\{i,j\}$, then their corresponding rows in $U$
	also agree on $S\cup S'$, by the definitions of $S$
	and $S'$.
	Again, they can differ only on $T$, a contradiction.
	
	Put
	\[
	p_1=\rank U[S],\quad\text{and}\quad
	p_2=\rank U[D\cup\{i,j\}].
	\]
	Since both matrices are reduced and have order at
	least two, we have $p_1,p_2\ge3$.
	All entries between $S\setminus\{i\}$ and
	$D\cup\{j\}$ are one.
	Taking the Schur complement at $i$ in
	$U[S\cup D\cup\{j\}]$ gives a block diagonal matrix.
	Consequently,
	\[
	 p_1+p_2-1=\rank U[S\cup D\cup\{j\}]\le k.
	\]
	In particular, $p_1,p_2\le k-2$, so the induction
	hypothesis gives
	\[
	a\le m(p_1-1)+1,\quad\text{and}\quad
	d+2\le m(p_2-1)+1.
	\]
	
	\medskip
	\noindent\textit{Case} 1: $p_1,p_2\ge4$.
	By \cite[Lemma 5(i)]{GMT} and the monotonicity of $m$,
	we have
	\[
	a+d+2
	\le m(p_1-1)+m(p_2-1)+2
	\le m(k-3)+3=B+2.
	\]
	This contradicts $a+d\ge B+1$.
	
	\medskip
	\noindent\textit{Case} 2: $p_2=3$.
	Then $d+2\le3$, giving $d=1$, as desired.
	
	\medskip
	\noindent\textit{Case} 3: $p_1=3$ and $p_2\ge4$.
	Since $p_2\le k-2$, we have $a\le3$ and $d+2\le B$.
	Together with $a+d\ge B+1$, this forces $a=3$,
	$d=B-2$, and $p_2=k-2$.
	Thus $U[D\cup\{i,j\}]$ attains equality in
	\Cref{thm:seidel} and is switching equivalent to
	$E_{k-2}$. Moreover, $D_{ij}$ contains all indices
	of this submatrix other than $i$ and $j$.
	If $k\ge7$, this contradicts
	\eqref{eq:no-full-difference}, since
	$E_{k-2}=\Phi(E_{k-4})$.
	
	Consequently, $k=6$, $B=4$, and $d=2$.
	The matrices $U[S]$ and $U[D\cup\{i,j\}]$ have
	ranks three and four, respectively, and both attain
	equality in \Cref{thm:seidel}.
	Since they are normalized at $i$, their binary blocks
	are $J_2-I_2$ and $J_3-I_3$, respectively.
	Hence $U[S\cup D\cup\{j\}]$ is the matrix $W$
	defined in \eqref{eq:exception}, up to simultaneous
	permutation.
	
	It remains to show that $N=8$.
	Since $N-q=a+d=5$ and $3=a\le q\le4$, we have
	$N\in\{8,9\}$.
	Suppose $N=9$. Then $q=4$ and $t=1$.
	Write $T=\{w\}$.
	
	The representative matrix $U[S\cup\{w\}]$ of $U[R]$
	has order four and rank at most four.
	By the induction hypothesis, its rank is exactly four,
	and its binary block is $J_3-I_3$.
	Therefore, when $w$ is added to the principal submatrix
	corresponding to $W$, its column in the resulting
	binary block has two ones on the first block of $C_0$
	and $c$ ones on the second block, where
	$c\in\{0,1,2\}$ because $u_{jw}=1$.
	
	Using \eqref{eq:J-I-inverse}, we find that the Schur complement of $C_0$ in the
	resulting binary matrix is
	$-(2+c^2/2-c)\ne0.$
	By \eqref{eq:normal-form}, it follows that
	$U[S\cup D\cup\{j,w\}]$ has rank seven,
	contradicting $k=6$.
	Thus $N=8$, completing the proof.
\end{proof}

\subsection{Completing the proof}
\label{sec:complete-rank}

\begin{proof}[\bfseries Proof of \Cref{thm:seidel}]
We may assume that $U$ is reduced by taking a
representative matrix.
Let $k=\rank U$ and let $N$ be the order of $U$.
We proceed by induction on $k$.

\medskip
\noindent\textit{Base cases}.
If $U$ has rank one, then its order is one. Moreover, $U$ cannot have rank two, as noted in \Cref{sec:prelim}. For $k=3,4$, \Cref{lem:binary-space:dimension} gives $s\le2$, so \eqref{eq:binary-bound} gives $N^2\le4k$, implying $N=k=m(k-1)+1$. In the normalization at any index, the binary block in \eqref{eq:normal-form} is a nonsingular symmetric zero-diagonal binary matrix of order two or three. It must be $J-I$ (for order three its determinant is twice the product of the three off-diagonal entries above the diagonal). Thus $U$ is switching equivalent to $E_k$.

\medskip
\noindent\textit{Order bound}.
Let $k\ge5$ and assume that the theorem holds for matrices of rank less than $k$. Put $B=m(k-3)+1$, so $m(k-1)=2B$. Since $m(k-1)+1\ge k$, any reduced shifted Seidel matrix of rank $k$ and order greater than $m(k-1)+1$ has a principal submatrix of rank $k$ and order $m(k-1)+2$. Since $U$ is reduced and this principal submatrix has the same rank as $U$, it is also reduced, as noted in \Cref{sec:prelim}. Thus it suffices to consider $N\in\{2B+1,2B+2\}$. By \Cref{thm:tau-cases}, $d=\delta_*(U)\le B-1$. Choose $D=D_{ij}$ of size $d$ and $R=Y\setminus D$, and normalize at~$i$. Then $N-d>B$, so \Cref{lem:minimum-deletion} applies, giving $q(U[R])\le B\le N-B-1$. By \Cref{lem:small-tau}, $d=1$ (the exceptional case does not occur). Consequently
\[
 N-1=|R|\le2q(U[R])\le2B.
\]
This proves the bound.

\medskip
\noindent\textit{Equality characterization}.
If $N=m(k-1)+1$, then all row classes of $U[R]$ have size two, and their representative matrix $P$ has order $B$ and rank $k-2$, by the induction hypothesis and $m(k-4)+1<B$. By \Cref{lem:minimum-deletion:classes}, ordering first the members that are one on $D$, then their partners, then the single index in $D$, gives exactly $\Phi(P)$. By induction $P$ is switching equivalent to $E_{k-2}$, so by \Cref{lem:phi:switching}, $U$ is switching equivalent to $E_k$. Conversely, by \Cref{lem:phi:rank,lem:phi:reduced}, $E_k$ is reduced of rank $k$ and order $m(k-1)+1$, and hence attains the bound.
\end{proof}


\begin{proof}[{\bfseries Proof of \Cref{conj:rank}}]
	Let $G$ be a reduced graph of rank $r\ge2$.
	Since $U(\widehat G)$ is reduced of rank $r+1$ and
	order $|V(G)|+1$, \Cref{thm:seidel} implies
	$|V(G)|\le m(r)$.
	
	It remains to characterize the graphs attaining equality.
	By \Cref{thm:seidel}, these are precisely the graphs
	obtained by normalizing $E_{r+1}$ at some index and then
	deleting that index. For $r\in\{2,3\}$, this gives only $K_r$.
	
	Suppose $r\ge4$ and put $P=E_{r-1}$, so that
	$E_{r+1}=\Phi(P)$. There are two possibilities,
	depending on the index at which $\Phi(P)$ is normalized.
	
	First, consider normalization at either copy of an index
	$i$ of $P$. Let $F$ be the extremal graph of rank $r-2$
	obtained by normalizing $P$ at $i$ and deleting $i$.
	Let $D_i$ be the diagonal matrix obtained from row $i$
	of $P$. Normalizing $\Phi(P)$ at either copy of $i$
	gives
	\[
	\begin{pmatrix}
	Q&Q&\mathbf v\\
	Q&Q&-\mathbf v\\
	\mathbf v^\top&-\mathbf v^\top&1
	\end{pmatrix},
	\]
	where $Q=D_iPD_i$ is $U(\widehat F)$ up to simultaneous
	permutation, and $\mathbf v=\pm D_i\mathbf1$.
	Interchanging the members of each pair for which
	$v_j=-1$ gives $\Phi(Q)$, with the normalization index
	in the first copy of $Q$. Deleting that index gives
	$U(F^{\rm a})$ up to simultaneous permutation, by
	\eqref{eq:constructions}.
	
	Next, normalization at the last index of $\Phi(P)$
	and deletion of that index gives $P^\pm$.
	Since $P$ is switching equivalent to $U(\widehat F)$
	for every extremal graph $F$ of rank $r-2$, we have
	$P^\pm=U(\widehat F)^\pm$ up to simultaneous permutation.
	By \eqref{eq:constructions}, the latter matrix is
	$U(F^{\rm c})$ up to simultaneous permutation.
	Thus Construction~{\rm(c)} gives a single isomorphism
	type, independently of the choice of $F$.
	
	Consequently, the extremal graphs of rank $r$ are
	precisely those obtained by applying Constructions
	{\rm(a)} and {\rm(c)} to the extremal graphs of rank
	$r-2$.
	
	Finally, Construction~{\rm(a)} is injective on isomorphism
	classes. Indeed, the new pendant vertex is unique.
	Deleting it and its neighbour leaves two copies
	of each vertex of $F$, and identifying vertices
	with identical neighbourhoods recovers $F$.
	Construction~{\rm(c)} gives a regular graph with no
	pendant vertices, whereas Construction~{\rm(a)} always
	produces a graph with a pendant vertex.
	Hence the number of extremal graphs up to isomorphism
	increases by exactly one when $r$ increases by two.
	Since there is a unique extremal graph for each of
	$r=2$ and $r=3$, this number is $\lfloor r/2\rfloor$.
\end{proof}

\section{Rank-one perturbations and proof of \Cref{conj:corank}}
\label{sec:corank}

Using \Cref{thm:seidel}, we first show that certain
integral extensions of shifted Seidel matrices
necessarily increase the rank. We then determine the maximum
order of a shifted Seidel matrix $U$ with distinct rows
in terms of the rank of $U+J$ and characterize the
matrices attaining equality. These results yield
\Cref{conj:corank}.

\subsection{Integral extensions}

We say that a shifted Seidel matrix $U$ admits an
\emph{odd extension} if there is an integral column
vector $\mathbf y$ with all entries odd and an integer
$c\equiv3\pmod4$ such that
\[
 \rank\begin{pmatrix}U&\mathbf y\\\mathbf y^\top&c\end{pmatrix}=\rank U.
\]
Existence of such an extension is invariant under switching, by applying the same sign congruence and permutation to the old indices, leaving the new index unchanged. Moreover, if $U$ admits an odd extension, then so does
every principal submatrix of $U$ having the same rank
as $U$.

\begin{lemma}\label{lem:integral-q}
A reduced shifted Seidel matrix of rank $k\ge3$ and order at least $m(k-1)$ admits no odd extension.
\end{lemma}
\begin{proof}
Let $N$ be the order of $U$ and $Y$ its index set.
We first rule out odd extensions for $E_3,E_4$ and $W$. After normalization at the first index, these matrices have the form \eqref{eq:normal-form} with $C=J_2-I_2$, $J_3-I_3$, or $C_0$ respectively. If $\mathbf y=(a,\mathbf z^\top)^\top$ is an odd integral column and $c\equiv3\pmod4$, the Schur complement of the first entry in the extension, divided by $-2$, is
\[
 \begin{pmatrix}C&\mathbf x\\\mathbf x^\top&e\end{pmatrix},
 \quad\text{where}\quad
 \mathbf x=(a\mathbf1-\mathbf z)/2,
 \quad\text{and}\quad e=(a^2-c)/2.
\]
Here $\mathbf x$ is integral and $e$ is odd. Rank preservation would imply $e=\mathbf x^\top C^{-1}\mathbf x$. For $C=J_l-I_l$, $l=2,3$, using
\eqref{eq:J-I-inverse}, we obtain
$e= \frac1{l-1}(\sum_i x_i)^2-\sum_i x_i^2.$
For $l=2$ it is even; for $l=3$ it is an even integer or a nonintegral half-integer, according to the parity of $\sum_i x_i$. For $C_0$ the two block contributions have the same properties. Hence none of these matrices admits an odd extension. This also proves the lemma for $k=3,4$ by \Cref{thm:seidel} and $N\ge k$.

We proceed by induction on $k$, the cases $k=3,4$
having been established above. Let $k\ge5$ and assume
that the lemma holds for all smaller ranks at least three.
Suppose $U$ admits an odd extension and put $B=m(k-3)+1$.
By \Cref{thm:seidel}, its order is $N\in\{2B,2B+1\}$.
Choose $D=D_{ij}$ attaining $\delta_*(U)=d$,
put $R=Y\setminus D$, and normalize at $i$.
By \Cref{thm:tau-cases}, $d\le B-1$, so $N-d>B$.

By \Cref{lem:minimum-deletion}, $q=q(U[R])\le B$.
If $q\le B-1$, then $q\le N-B-1$;
\Cref{lem:small-tau} yields $d=1$ or $k=6$ with
a principal submatrix switching equivalent to $W$.
The former is impossible since $N-1\le2q\le2B-2$,
and the latter is impossible since $W$ has full rank
six and admits no odd extension.

Thus $q=B$ and a representative matrix $P$ of $U[R]$
has rank $k-2$, by \Cref{thm:seidel} and
$m(k-4)+1<B$.
Choose its index set $L$ containing $i$ (and hence
not $j$), and choose $w\in D$.

Restrict an odd extension of $U$ to $L\cup\{j,w\}$
and the new index.
Subtracting row and column $i$ from row and column $j$,
and also from row and column $w$, gives a matrix of
rank at most $k$ of the form
\[
\begin{pmatrix}
P & \mathbf0 & \mathbf b & \mathbf y\\
\mathbf0^\top & 0 & -2 & a\\
\mathbf b^\top & -2 & 0 & b\\
\mathbf y^\top & a & b & c
\end{pmatrix},
\]
where $a$, $b$, and all entries of $\mathbf b$ are even,
$\mathbf y$ is odd, and $c\equiv3\pmod4$.
Taking the Schur complement of the block on $j,w$ gives
\[
\begin{pmatrix}
P & \mathbf y+\frac{a}2\mathbf b\\
\mathbf y^\top\!+\frac{a}2\mathbf b^\top & c+ab
\end{pmatrix}.
\]
This matrix has rank at most $k-2$, and hence exactly
$k-2$, since $\rank P=k-2$.
Thus it is an odd extension of $P$.
Since $P$ is reduced and has order $m(k-3)+1$,
this contradicts the induction hypothesis.
\end{proof}

\subsection{The perturbation bound}

\begin{theorem}\label{thm:perturbation}
	Let $U$ be a shifted Seidel matrix with distinct rows
	and $\rank(U+J)=r\ge2$. Then its order is at most
	$m(r)$. For $r\ge4$, equality holds if and only if
	$U$ is equal to $E_{r-1}^\pm$ up to simultaneous
	permutation of rows and columns.
\end{theorem}

\begin{proof}
Let $k=\rank U$ and let $N$ be the order of $U$.
Since $\rank(U+J)=r$, we have $k\le r+1$.

\medskip
\noindent\textit{Case} 1: Two rows of $U$ differ only in sign.
Their corresponding rows in $U+J$ sum to
$2\mathbf1^\top$. Hence the row space of $U+J$
is spanned by the rows of $U$ and $\mathbf1^\top$.
Since the two corresponding columns of $U$ also
differ only in sign, every vector in the row space
of $U$ has opposite entries at these two coordinates.
Thus $\mathbf1^\top$ does not belong to the row
space of $U$, and hence $k=r-1$.

Since the rows of $U$ are distinct, every row class
has at most two members.
If $r=2$, then $k=1$ and $q(U)=1$, giving
$N\le2=m(2)$.
The case $r=3$ is impossible, since $U$ cannot
have rank two.
For $r\ge4$, by \Cref{thm:seidel},
\[
 N\le2q(U)\le2(m(r-2)+1)=m(r).
\]
If equality holds, then every row class has two
members and, by \Cref{thm:seidel}, a representative
matrix $P$ of $U$ is switching equivalent to $E_{r-1}$.
Thus $U$ is equal to $P^\pm$, and hence to
$E_{r-1}^\pm$, up to simultaneous permutation,
as noted in \Cref{sec:prelim}.
Conversely, $E_{r-1}^\pm$ has distinct rows,
rank $r-1$, and order $m(r)$.
Since it contains two rows that differ only in sign,
the rank of its sum with $J$ is $r$, as shown above.

\medskip
\noindent\textit{Case} 2: $U$ is reduced.
If $k=1$, then $N=1$ and the assertion is immediate.
Thus $k\ge3$.
By \Cref{thm:seidel}, $N\le m(r)+1$.
Suppose $N=m(r)+1$, or $r\ge4$ and $N=m(r)$.
Since $N>m(r-1)+1$ in either case,
\Cref{thm:seidel} implies $k=r+1$.
But
\[
 \rank\begin{pmatrix}
 U&\mathbf1\\
 \mathbf1^\top&-1
 \end{pmatrix}
 =1+\rank(U+J)=r+1=k,
\]
giving an odd extension of $U$, contrary to
\Cref{lem:integral-q}.
\end{proof}


\begin{proof}[\bfseries Proof of \Cref{conj:corank}]
Let $G$ be a coreduced graph of corank $r\ge2$.
Since
\[
 U(\overline G)+J=2(A(G)+I),
\]
the matrix $U(\overline G)$ has distinct rows and
$\rank(U(\overline G)+J)=r$.
Hence \Cref{thm:perturbation} implies
$|V(G)|\le m(r)$.

For $r\ge4$, equality holds if and only if
$U(\overline G)$ is equal to $E_{r-1}^\pm$ up to
simultaneous permutation of rows and columns.
As shown in the proof of \Cref{conj:rank}, this is
equivalent to $\overline G$ being the unique regular extremal graph
of rank $r$, obtained from $K_2$ or $K_3$ by repeated
applications of Construction~{\rm(c)}.

For $r=2,3$, the equality cases follow by inspecting
the graphs on two and three vertices, respectively.
\end{proof}

\paragraph{AI statement.}
The ideas of using Seidel matrices and incorporating
binary rank arguments originated with the author.
ChatGPT (GPT-6.1 and GPT-6.2) assisted in developing these ideas.
The author wrote the proofs and used ChatGPT to
polish the exposition.


\begin{thebibliography}{99}
\bibitem{ACK} S. Akbari, P. J. Cameron and G. B. Khosrovshahi, \emph{Ranks and signatures of adjacency matrices}, unpublished preprint, (2004), \url{https://maths.qmul.ac.uk/~pjc/preprints/ranksign.pdf}.

\bibitem{AS} N. Alon and P. D. Seymour, A counterexample to the rank-coloring conjecture, \emph{J. Graph Theory} \textbf{13} (1989), 523--525.

\bibitem{BLS} G. Beniamini, N. Linial and A. Shraibman, The Rank-Ramsey problem and the Log-Rank conjecture, \emph{Combin. Probab. Comput.} \textbf{35} (2026), 568--597.

\bibitem{BR} F. K. Bell and P. Rowlinson, On the multiplicities of graph eigenvalues, \emph{Bull. London Math. Soc.} \textbf{35} (2003), 401--408.

\bibitem{EGcorank} H. Esmailian and E. Ghorbani, Maximum order of graphs with a given corank, \emph{Linear Algebra Appl.} \textbf{588} (2020), 122--133.

\bibitem{EGHK} H. Esmailian, E. Ghorbani, S. Hossein Ghorban and G. B. Khosrovshahi, Maximal graphs with respect to rank, \emph{Discrete Math.} \textbf{344} (2021), 112191.

\bibitem{EGLO} H. Esmailian and E. Ghorbani, A discrete variation of the Littlewood--Offord problem, \emph{Electron. J. Combin.} \textbf{31} (2024), \#P1.42.

\bibitem{GMT} E. Ghorbani, A. Mohammadian and B. Tayfeh-Rezaie, On order and rank of graphs, \emph{Combinatorica} \textbf{35} (2015), 655--668.

\bibitem{GMTbip} E. Ghorbani, A. Mohammadian and B. Tayfeh-Rezaie, Maximum order of trees and bipartite graphs with a given rank, \emph{Discrete Math.} \textbf{312} (2012), 3498--3501.

\bibitem{GMTtri} E. Ghorbani, A. Mohammadian and B. Tayfeh-Rezaie, Maximum order of triangle-free graphs with a given rank, \emph{J. Graph Theory} \textbf{79} (2015), 145--158.

\bibitem{GPW} M. G\"o\"os, T. Pitassi and T. Watson, Deterministic communication vs. partition number, \emph{SIAM J. Comput.} \textbf{47} (2018), 2435--2450.

\bibitem{HJ} R. A. Horn and C. R. Johnson, \emph{Matrix Analysis}, 2nd ed., Cambridge University Press, Cambridge, 2013.

\bibitem{HP} W. H. Haemers and M. J. P. Peeters, The maximum order of adjacency matrices of graphs with a given rank, \emph{Designs, Codes and Cryptography} \textbf{65} (2012), 223--232.

\bibitem{KL} A. Kotlov and L. Lov\'asz, The rank and size of graphs, \emph{J. Graph Theory} \textbf{23} (1996), 185--189.

\bibitem{Lovett} S. Lovett, Communication is bounded by root of rank, \emph{J. ACM} \textbf{63} (2016), Article 1.

\bibitem{LS} L. Lov\'asz and M. Saks, Lattices, M\"obius functions and communication complexity, in: \emph{Proceedings of the 29th Annual IEEE Symposium on Foundations of Computer Science}, 1988, pp. 81--90.

\bibitem{LS2} L. Lov\'asz and M. Saks, Communication complexity and combinatorial lattice theory, \emph{J. Comput. System Sci.} \textbf{47} (1993), 322--349.

\bibitem{NW} N. Nisan and A. Wigderson, On rank vs. communication complexity, \emph{Combinatorica} \textbf{15} (1995), 557--565.

\bibitem{Raz} A. A. Razborov, The gap between the chromatic number of a graph and the rank of its adjacency matrix is superlinear, \emph{Discrete Math.} \textbf{108} (1992), 393--396.

\bibitem{Row} P. Rowlinson, On graphs with multiple eigenvalues, \emph{Linear Algebra Appl.} \textbf{283} (1998), 75--85.

\bibitem{ST} B. Sudakov and I. Tomon, Matrix discrepancy and the log-rank conjecture, \emph{Math. Program.} \textbf{212} (2025), 567--579.

\bibitem{TR} B. Tayfeh-Rezaie, \emph{Lecture Notes on the Star Complement Technique}, 2009, available at \url{http://math.ipm.ac.ir/~tayfeh-r/papersandpreprints/starcompl.pdf}.

\bibitem{VN} C. van Nuffelen, A bound for the chromatic number of a graph, \emph{Amer. Math. Monthly} \textbf{83} (1976), 265--266.
\end{thebibliography}
\end{document}